\documentclass[12pt]{article}

 \usepackage{amscd,amsfonts,amssymb,amsthm,amscd,amsmath,latexsym}
\begin{titlepage}

\title{Determinant bundles, boundaries, and surgery}
\author{Ulrich Bunke and Jinsung Park}
\end{titlepage}
\begin{document}
\maketitle
\begin{abstract}
In this note we specialize and illustrate the ideas developed
in the paper \cite{bunke020} {\em Families, Eta forms, and Deligne cohomology}
in the case of the determinant line bundle.
We discuss the surgery formula in the adiabatic limit  
using the adiabatic decomposition formula
of the zeta regularized determinant of the Dirac Laplacian in \cite{PW2}. 
\end{abstract}

\maketitle
\newcommand{\Fred}{{\tt Fred}}
\newcommand{\Comp}{{\tt K}}
\newcommand{\vC}{\v{C}}

\newcommand{\DD}{\mathbf{v}}
\newcommand{\Z}{\mathbb{Z}}
\newcommand{\orient}{{\rm or}}
\newcommand{\cSet}{\mathcal{S}et}
\newcommand{\by}{{\bf y}}
\newcommand{\bE}{{\bf E}}
\newcommand{\Face}{{\tt Face}}
\newcommand{\cDelta}{\mathbf{\Delta}}
\newcommand{\LIM}{{\tt LIM}}
\newcommand{\diag}{{\tt diag}}
 \newcommand{\dist}{{\tt dist}}
\newcommand{\kaaa}{{\frak k}}
\newcommand{\paaa}{{\frak p}}
\newcommand{\vp}{{\varphi}}
\newcommand{\taaa}{{\frak t}}
\newcommand{\haaa}{{\frak h}}
\newcommand{\R}{\mathbb{R}}
\newcommand{\Hh}{{\bf H}}
\newcommand{\Rep}{{\tt Rep}}
\newcommand{\Hb}{mathbb{H}}
\newcommand{\Q}{\mathbb{Q}}
\newcommand{\str}{{\tt str}}
\newcommand{\Ind}{{\tt ind}}
\newcommand{\triv}{{\tt triv}}
\newcommand{\bD}{{\bf D}}
\newcommand{\bF}{{\bf F}}
\newcommand{\tX}{{\tt X}}
\newcommand{\Cliff}{{\tt Cliff}}
\newcommand{\tY}{{\tt Y}}
\newcommand{\tZ}{{\tt Z}}
\newcommand{\tV}{{\tt V}}
\newcommand{\tR}{{\tt R}}
\newcommand{\bu}{{\mathbf{u}}}
\newcommand{\Fam}{{\tt Fam}}
\newcommand{\Cusp}{{\tt Cusp}}
\newcommand{\bT}{{\bf T}}
\newcommand{\bK}{{\bf K}}
\newcommand{\bo}{{\bf o}}
\newcommand{\K}{\mathbb{K}}
\newcommand{\tH}{{\tt H}}
\newcommand{\bS}{\mathbf{S}}
\newcommand{\bB}{\mathbf{B}}
\newcommand{\tW}{{\tt W}}
\newcommand{\tF}{{\tt F}}
\newcommand{\bA}{\mathbf{A}}
\newcommand{\bL}{{\bf L}}
 \newcommand{\bom}{\mathbf{\Omega}}
\newcommand{\bundle}{{\tt bundle}}
\newcommand{\ch}{\mathbf{ch}}
\newcommand{\ve}{{\varepsilon}}
\newcommand{\C}{\mathbb{C}}
\newcommand{\gen}{{\tt gen}}
\newcommand{\cTop}{\mathcal{T}op}
\newcommand{\bP}{\mathbf{P}}
\newcommand{\Naaa}{\mathbf{N}}
\newcommand{\image}{{\tt image}}
\newcommand{\gaaa}{{\frak g}}
\newcommand{\zaaa}{{\frak z}}
\newcommand{\saaa}{{\frak s}}
\newcommand{\laaa}{{\frak l}}
\newcommand{\bN}{\mathbf{N}}
\newcommand{\stimes}{{\times\hspace{-1mm}\bf |}}
\newcommand{\ausg}{{\rm end}}
\newcommand{\bff}{{\bf f}}
\newcommand{\maaa}{{\frak m}}
\newcommand{\aaaa}{{\frak a}}
\newcommand{\naaa}{{\frak n}}
\newcommand{\brr}{{\bf r}}
\newcommand{\res}{{\tt res}}
\newcommand{\Aut}{{\tt Aut}}
\newcommand{\Pol}{{\tt Pol}}
\newcommand{\Tr}{{\tt Tr}}
\newcommand{\cT}{\mathcal{T}}
\newcommand{\dom}{{\tt dom}}
\newcommand{\Line}{{\tt Line}}
\newcommand{\db}{{\bar{\partial}}}
\newcommand{\Sf}{{\tt  Sf}}
\newcommand{\g}{{\gaaa}}
\newcommand{\Det}{{\tt Det}}
\newcommand{\cZ}{\mathcal{Z}}
\newcommand{\cH}{\mathcal{H}}
\newcommand{\cM}{\mathcal{M}}
\newcommand{\interi}{{\ttint}}
\newcommand{\singsupp}{{\tt singsupp}}
\newcommand{\cE}{\mathcal{E}}
\newcommand{\ccR}{\mathcal{R}}
\newcommand{\hol}{{\tt hol}}
\newcommand{\cV}{\mathcal{V}}
\newcommand{\cY}{\mathcal{Y}}
\newcommand{\cW}{\mathcal{W}}
\newcommand{\dR}{{\tt dR}}
\newcommand{\del}{{\tt del}}
\newcommand{\bdel}{\mathbf{del}}
\newcommand{\cI}{\mathcal{I}}
\newcommand{\cC}{\mathcal{C}}
\newcommand{\cK}{\mathcal{K}}
\newcommand{\cA}{\mathcal{A}}
\newcommand{\cU}{\mathcal{U}}
\newcommand{\Hom}{{\tt Hom}}
\newcommand{\vol}{{\tt vol}}
\newcommand{\cO}{\mathcal{O}}
\newcommand{\End}{{\tt End}}
\newcommand{\Ext}{{\tt Ext}}
\newcommand{\rk}{{\tt rank}}
\newcommand{\im}{{\tt im}}
\newcommand{\sign}{{\tt sign}}
\newcommand{\spann}{{\tt span}}
\newcommand{\symm}{{\tt symm}}
\newcommand{\cF}{\mathcal{F}}
\newcommand{\cD}{\mathcal{D}}
\newcommand{\bC}{\mathbf{C}}
\newcommand{\bbeta}{\mathbf{\eta}}
\newcommand{\bOmega}{\mathbf{\Omega}}
\newcommand{\bbz}{{\bf z}}
\newcommand{\bc}{\mathbf{c}}
\newcommand{\bb}{\mathbf{b}}
\newcommand{\bd}{\mathbf{d}}
\newcommand{\Ree}{{\tt Re }}
\newcommand{\tU}{{\tt U}}
\newcommand{\Res}{{\tt Res}}
\newcommand{\Imm}{{\tt Im}}
\newcommand{\inter}{{\tt int}}
\newcommand{\clo}{{\tt clo}}
\newcommand{\tg}{{\tt tg}}
\newcommand{\ee}{{\tt e}}
\newcommand{\Li}{{\tt Li}}
\newcommand{\cN}{\mathcal{N}}
 \newcommand{\conv}{{\tt conv}}
\newcommand{\op}{{\tt Op}}
\newcommand{\tr}{{\tt tr}}
\newcommand{\ctg}{{\tt ctg}}
\newcommand{\degg}{{\tt deg}}
\newcommand{\Ad}{{\tt  Ad}}
\newcommand{\ad}{{\tt ad}}
\newcommand{\codim}{{\tt codim}}
\newcommand{\Gr}{{\tt Gr}}
\newcommand{\coker}{{\tt coker}}
\newcommand{\id}{{\tt id}}
\newcommand{\tL}{{\tt L}}
\newcommand{\ord}{{\tt ord}}
\newcommand{\nat}{\mathbb{N}}
\newcommand{\supp}{{\tt supp}}
\newcommand{\sing}{{\tt sing}}
\newcommand{\spec}{{\tt spec}}
\newcommand{\Ann}{{\tt Ann}}
 \newcommand{\Or}{{\tt Or }}
\newcommand{\Diff}{\mathcal{D}iff}
\newcommand{\cB}{\mathcal{B}}
\def\imath{{i}}
\newcommand{\cR}{\mathcal{R}}
\def\hB{\hspace*{\fill}$\Box$ \newline\noindent}
\newcommand{\varho}{\varrho}
\newcommand{\ind}{{\tt index}}
\newcommand{\Indu}{{\tt Ind}}
\newcommand{\Fin}{{\tt Fin}}
\newcommand{\cS}{\mathcal{S}}
\newcommand{\orig}{\mathcal{O}}
\def\hB{\hspace*{\fill}$\Box$ \\[0.5cm]\noindent}
\newcommand{\cL}{\mathcal{L}}
 \newcommand{\cG}{\mathcal{G}}
\newcommand{\Mat}{{\TT Mat}}
\newcommand{\cP}{\mathcal{P}}
\newcommand{\bv}{\mathbf{v}}
\newcommand{\cQ}{\mathcal{Q}}
 \newcommand{\cX}{\mathcal{X}}
\newcommand{\bH}{\mathbf{H}}
\newcommand{\bW}{\mathbf{W}}
\newcommand{\pr}{{\tt pr}}
\newcommand{\bX}{\mathbf{X}}
\newcommand{\bY}{\mathbf{Y}}
\newcommand{\bZ}{\mathbf{Z}}
\newcommand{\bz}{\mathbf{z}}
\newcommand{\bkappa}{\mathbf{\kappa}}
\newcommand{\ev}{{\tt ev}}
\newcommand{\bV}{\mathbf{V}}
\newcommand{\Gerbe}{{\tt Gerbe}}
\newcommand{\gerbe}{{\tt gerbe}}
\newcommand{\hA}{\mathbf{\hat A}}
\newcommand{\Sets}{\mathcal{S}ets}
\newcommand{\bg}{{\mathbf{g}}}
\def\det{{\tt det}}

\newtheorem{prop}{Proposition}[section]
\newtheorem{lem}[prop]{Lemma}
\newtheorem{ddd}[prop]{Definition}
\newtheorem{theorem}[prop]{Theorem}
\newtheorem{kor}[prop]{Corollary}
\newtheorem{ass}[prop]{Assumption}
\newtheorem{con}[prop]{Conjecture}
\newtheorem{prob}[prop]{Problem}
\newtheorem{fact}[prop]{Fact}

\tableofcontents

\parskip3ex

\section{Introduction}

Recently, the hirachy of geometric objects over a smooth manifold which starts with $U(1)$-valued smooth functions, hermitean line bundles with connections, and geometric gerbes in degree zero, one, and two, respectively, became a subject of intensive investigation. One of the motivation stems from the philosophy  that an object of degree $k+1$ on a base gives rise, via transgression, to an object of degree $k$ on the free loop space of the base. So one wants to study low degree objects like line bundles on iterated loop spaces via higher degree objects on the original much smaller base.

From a categorial point of view  the objects in degree zero form just a set. In degree one we have a category of objects. A proper study of gerbes already requires
(weak) $2$-categories, and the higher objects in this hirachy require higher categories. Currently there are many different notions of higher categories (see e.g. \cite{leinster01}), and the nature of the higher geometric objects has to be clarified.

The different levels in this hirachy should be linked as follows : There is a monoidal structure (a sort of tensor product) on the objects in degree $k$ such that each object is invertible. 
The geometric objects in degree $k$
act ``transitively'' on the morphisms between objects in degree $k+1$ such that the monoidal structure relates to composition.
 For example,  an $U(1)$-valued function acts on the space of morphisms between hermitean line bundles.

Already for gerbes we have several geometric pictures, e.g.  as sheafs of categories (Brylinski \cite{brylinski93}),  as bundle gerbes (Murray \cite{murray96}), or as higher line bundles (Hitchin \cite{hitchin99}). A geometric theory for the higher degree has been proposed in (Gajer \cite{gajer97}).

Though there may be different geometric realizations there is no doubt that the
isomorphism classes of the objects in degree $k$ over the base $B$ are classified by the smooth Deligne cohomology $H^k_{Del}(B)$. This cohomology group is defined as the $k-1$`th  \v{C}ech hypercohomology group of the complex of sheaves
$$\cR^{k-1}_B : 0\rightarrow \underline{U(1)}_B^\infty\stackrel{\frac{1}{2\pi i}d\log }{\longrightarrow}\cA^1_B\rightarrow\dots\rightarrow \cA^{k-1}_B\rightarrow 0\ .$$
Here $\underline{U(1)}_B^\infty$ denotes the sheaf of smooth $U(1)$-valued functions, and $\cA_B^p$ denote the sheaves of real smooth differential forms. 

The \v{C}ech complex provides a very simple geometric picture of the hirachy where we can completely avoid the language higher categories. An object in degree $k+1$ is just a \v{C}ech $k$-cocycle $\bc$ for $\cR^k_B$. The tensor product of objects  is the sum of cocycles.
A morphism $\bb:\bc_0\rightarrow \bc_1$ is just given by a 
\v{C}ech $k-1$-cochain $\bb$ such that $\bd\bb=\bc_1-\bc_0$.
The difference between two such morphisms $\bb_0-\bb_1$ can be considered as
a \v{C}ech $k-1$-cocycle for $\cR^{k-1}_B$ in a natural way.
Note that all morphisms are invertible, and that the isomorphism classes
of objects in degree $k+1$ are indeed classified by $H^{k+1}_{Del}(B)$.
A formula for the transgression was given by Gawedski \cite{brylinski93}.
A unit norm section of an object $\bc$ of degree $k+1$ is by definition a $k$-chain $\bb$
such that $\bd\bb+\bc$ is concentrated in \v{C}ech degree zero.
This sum comes in fact from a global $k$-form which will also be denoted by
$\nabla^{\bc}\bb$. The difference of two sections is a degree $k$-object in a natural way. A morphism from the zero object to $\bc$ is the same  as a ``parallel'' section.

By a refinement of local index theory a geometric family of Dirac operators 
(We refer to  \cite{bunke020} for a definition of this notion. It essentially includes all information needed for doing local index theory.)
may give rise to objects in this hirachy. E.g. in degree zero we have the exponentiated $\eta$-invariant of a family of Dirac operators on odd-dimensional manifolds. The most prominent example is the determinant line bundle 
which is associated to a geometric family of Dirac operators on even dimensional manifolds. This line bundle  is equipped with the Quillen metric and the Bismut-Freed connection (see \cite{berlinegetzlervergne92} for details).
Recently Lott constructed an isomorphism class of a geometric index gerbe which is again associated to a family of Dirac operators on odd-dimensional manifolds. Under certain conditions one can also construct isomorphism classes of  higher degree objects
which however depend on additional choices (see e.g. Lott \cite{lott01} and \cite{bunke020}).
 
The isomorphism class of the index gerbe of Lott was first constructed using Hitchin's picture of gerbes. The geometric object depend on choices, and it is ony the isomorphism class which is independent of choices. Note in constrast, that in the case of the classical construction of the determinant bundle we really know a canonical realization of the fibre over any point of the base (see (\ref{er})). This makes a huge difference if we want to speak about sections or trivializations etc.

The higher degree objects in \cite{lott01}, and the objects in  all degrees in  \cite{bunke020} are constructed as \v{C}ech cocycles again depending on choices. In \cite{bunke020} these choices are called tamed resolutions.

The goal of the present paper is to demonstrate
the flexibility of the \v{C}ech complex picture.
We start with degree $k=1$.
So we provide a description in this language of classical constructions related
to the determinant bundle like the determinant section, the $\tau$-section associated to a zero bordism due to Dai-Freed \cite{daifreed94} and a surgery formula.  

Then we turn to the next degree $k=2$. We show that the \v{C}ech complex picture is strong enough view a generalization of the determinant line bundle on a family with boundary as a trivialization of the index gerbe associated to the boundary, and to provide a glueing formula for these generalizations.
This generalizes the Dai-Freed $\tau$-section associated to a zero bordism above to  degree $2$.

It is now straight forward to generalize this picture to all degrees.
The \v{C}ech complex picture is of course very simple and does not reflect
the true geometry nature of the higher objects. E.g. it is not clear that we have an integration over the fibre on the \v{C}ech complex level for fibres of dimension $\ge 2$, though such a transgression  exists on the level of isomorphism classes. It is also not so obvious
which picture will be eventually the natural one adapted to local index theory. In a truely higher categorial world trivializations
are a very complicated object to study.

\section{Hermitian line bundles with connection and Deligne cohomology}\label{lb}

Let $B$ be a smooth manifold. If $\bL=(L,h^L,\nabla^L)$ is a hermitian line bundle with connection on $B$, then by $\Gamma(L)$ and $\Gamma_1(\bL)$ we denote its spaces of sections and of unit norm sections. The set of isomorphism classes of hermitean line bundles forms an abelian group $\Line(B)$ where the group operation is induced by the tensor product.

Consider for $i=0,1$ hermitian line bundles with connection $\bL_i$ and sections $\phi_i\in\Gamma_1(\bL_i)$.
We say that the pairs $(\phi_i,\bL_i)$ are isomorphic iff there is an isomorphism $f:\bL_0\rightarrow \bL_1$ such that $f(\phi_0)=\phi_1$.
We define a product $(\phi_0,\bL_0)\otimes (\phi_1,\bL_1):=(\phi_0\otimes\phi_1,\bL_0\otimes \bL_1)$. The set of isomorphism classes of such pairs forms the abelian group $\Gamma\Line(B)$.  

Let $\cU:=(U_l)_{l\in I}$ be an open covering and $s:=(s_l)$, $s_l\in \Gamma_1(\bL_{|U_l})$ be a family of local sections. 
If $\cU^\prime:=(U^\prime_k)_{k\in J}$, $r:J\rightarrow I$, is a refinement, then
we obtain a family $s^\prime:=(s_k^\prime)$, $s_k^\prime:=s_{r(k)|U^\prime_k}$.
We say that the tuples $(\bL_0,\cU_0,s_0)$, $(\bL_1,\cU_1,s_1)$ are isomorphic, if there exists an isomorphism $f:\bL_0\rightarrow \bL_1$ such that
after a common refinement of the coverings $f(s_0)=s_1$.
We have again a tensor product of locally trivialized bundles.
By $\widetilde{\Line(B)}$ we denote the abelian group of isomorphism classes of
tuples $(\bL,\cU,s)$.

Finally, we introduce the abelian group $\widetilde{\Gamma \Line(B)}$ of isomorphism classes of tuples $(\bL,\phi,\cU,s)$, where $\phi\in\Gamma_1(\bL)$.
 
\underline{Remark :}
Note that the forgetful maps $\widetilde{\Line(B)}\rightarrow \Line(B)$
and $\widetilde{\Gamma \Line(B)}\rightarrow \Gamma\Line(B)$ are isomorphisms.

We consider the complex of sheaves 
$$\cR^1_B: 0\rightarrow \underline{U(1)}^\infty_B\stackrel{\frac{1}{2\pi i}d \log }{\rightarrow}\cA^1_B\rightarrow 0\ .$$
The degree-two Deligne cohomology $H^2_{Del}(B)$ is by definition the degree-one \v{C}ech hypercohomology of $\cR^1_B$. A cochain is given by a pair
$\bc=(c^{0,1},c^{1,0})\in \check{\bC}^1(B,\cR^1_B)$, where $c^{0,1}\in \check{C}^0(B,\cA^1_B)$ and
$c^{1,0}\in\check{C}^1(B,\underline{U(1)}^\infty_B)$.
The chain $\bc$ is closed if $\delta c^{0,1}=-\frac{1}{2\pi i} d\log c^{1,0}$.

We define a homomorphism $\widetilde{l_B}:\widetilde{\Line(B)}\rightarrow \check{\bZ}^1(B,\cR^1_B)$ by
$\widetilde{l_B}(\bL,\cU,s)=:\bc$ with
$c^{0,1}_l=\frac{1}{2\pi i} \nabla^L \log s_l$ and
$c^{1,0}_{lk}= \frac{s_l}{s_k}$.
It is easy to check that $\widetilde{l_B}$ is an isomorphism of groups.
It induces an isomorphism
$l_B:\Line(B)\rightarrow H^2_{Del}(B)$.

\underline{Remark :}
The remark above has a counterpart for degree-one \v{C}ech cocycles.
In fact, the projection
$$\check{\bZ}^1(B,\cR^1_B)\rightarrow H^2_{Del}(B)$$
is an isomorphism.

Let $\bc\in \check{\bZ}^1(B,\cR^1_B)$. A unit norm section of $\bc$ is by definition a chain $\bu\in\check{\bC}^0(B,\cR^1_B)$, $\bu=(u^{0,0})$, 
$u:=u^{0,0}\in  \check{C}^0(B, \underline{U(1)}^\infty_B)$ such that
$\delta u= c^{1,0}$. By $\Gamma_1(\bc)$ we denote the set of all such sections. We form the group
$\Gamma  \check{\bZ}^1(B,\cR^1_B)$ which consists of all pairs $(\bc,\bu)$, where $\bu\in\Gamma_1(\bc)$. Then we can extend $\widetilde{l_B}$ to an isomorphism
$\widetilde{\Gamma l_B}:\widetilde{\Gamma \Line(B)}\rightarrow \Gamma \check{\bZ}^1(B,\cR^1_B)$ such that
$\widetilde{\Gamma l_B}(\bL,\phi,\cU,s)=(\bc,\bu)$ with
$\bc=\widetilde{l_B}(\bL,\cU,s)$ and $u_l=\frac{\phi_{|U_l}}{s_l}$.

There are natural homomorphisms
$R:H^2_{Del}(B)\rightarrow \cA_B^2(B,d=0)$ (the curvature form) and
$\bv:H^2_{Del}(B)\rightarrow \check{H}^2(B,\underline{\Z}_B)$
(the first Chern class), which are given in terms of a representative $\bc$  by $R^{\bc}_{|U_l}=d c^{0,1}_l$ and
$\bv(\bc)=[\delta \frac{1}{2\pi i}\log c^{1,0}]$.

Note that the composition $\bv\circ l_B$ is equal  to the first Chern class $c_1:\Line(B)\rightarrow\check{H}^2(B,\underline{\Z}_B)$, while $R\circ l_B(\bL)=\frac{1}{2\pi i}R^{\nabla^L}$. 

Let $(\bL,\cU,s)$ be given. A section $\phi\in\Gamma(L)$ gives rise to
a family $(\phi_l)_{l\in I}$ of complex valued functions such that
$\phi_ls_l=\phi_{|U_l}$ for all $l$.
If $\bc=\widetilde{l_B}(\bL,\cU,s)$, then
any family $(\phi_l)$ such that
$\phi_k=c^{0,1}_{lk}\phi_l$ on $U_l\cap U_k$ for all $l,k\in I$
defines a section of $L$.

The one-form $\nabla^L \log \phi$ is given in terms of the $\phi_l$ by
$\nabla^L \log \phi_{|U_l}=d  \log \phi_l + 2\pi i c^{0,1}_l$.

In Section \ref{bound} we introduce a similar translation between gerbes and  degree-three Deligne cohomology.

\section{The determinant line bundle}

We consider a geometric family $\cE$ (see \cite{bunke020}) with closed even-dimensional fibers over a manifold $B$. This notion combines the data of a smooth fiber bundle $\pi:E\rightarrow B$ with closed even-dimensional fibers which is equipped with an orientation of the vertical bundle $T^v\pi$, a vertical metric $g^{T^v\pi}$, a horizontal distribution $T^h\pi$, and a family of $\Z_2$-graded Dirac bundles $\cV:=(V,h^V,\nabla^V,c)$.

Let $D^\pm(\cE)=(D^\pm_b)_{b\in B}$ be the associated family of chiral Dirac operators. We assume that $\ind(D^+_b)=0$ for all $b\in B$.
The determinant line bundle of $D^+(\cE)$ is a hermitian line bundle with connection over $B$ which will be denoted by $\det(\cE)$. The fiber of $\det(\cE)$ over $b\in B$ is canonically isomorphic to
\begin{equation}\label{er}\det(\cE)_b\cong \Hom(\Lambda^{max}\ker(D^+_b), \Lambda^{max}\ker(D^-_b)) \ .\end{equation} The determinant line bundle is equipped with the Quillen metric
and the Bismut-Freed connection. We refer to \cite{berlinegetzlervergne92}, Ch.9,  for the construction of the metric and the connection.

The bundle $\det(\cE)$ comes with a canonical section $\Det(\cE)$. It is given by
$$\Det(\cE)(b):=\left\{ \begin{array}{cc}1& \ker(D^+_b)=0\\
0 & \ker(D^+_b)\not=0\end{array}\right.\ .$$
Note that this definition makes sense because
$\det(\cE)_b\cong \C$ canonically if $\ker(D^+_b)=0$.

Let $K^0_2(B)\subset K^0(B)$ be the subgroup which consists of all classes which vanish on the one-skeleton of $B$.
Since we assume that the index of the chiral Dirac operator on each fiber
vanishes we have $\ind(\cE)\in K_2^0(B)$, where
$\ind(\cE)\in K^0(B)$ denotes the index of the family $D^+(\cE)$ 

In \cite{bunke020} we constructed a refined index $\ind^2_{Del}(\cE)\in H^2_{Del}(B)$
such that
$$R^{\ind^2_{Del}(\cE)}=\Omega^2(\cE)$$
and $$\bv(\ind^2_{Del}(\cE)) = c_1(\ind(\cE))\ .$$
Here $\Omega^2(\cE)$ denotes the degree-two component of the local index form which is given by 
$$\Omega(\cE):=\int_{E/B}\hA(\nabla^{T^v\pi})\ch(\nabla^{\bW})\ ,$$
where $\bW$ is the (locally defined) twisting bundle with induced connection, and $\nabla^{T^v\pi}$ is the connection induced by the vertical metric and the horizontal distribution. Furthermore we have verified 
by holonomy comparison that
$$l_B(\det(\cE))=\ind^2_{Del}(\cE)\ .$$

\section{Local trivializations of the determinant line bundle}\label{loctr}

In this section we want to lift the correspondence 
$$l_B(\det(\cE))=\ind^2_{Del}(\cE)$$ to the level of locally trivialized bundles and cycles. To this end we recall the construction of cycles $\bc\in  \check{\bZ}^1(B,\cR^1_B)$ which represent $\ind^2_{Del}(\cE)$, and then we exhibit
the corresponding local trivialization $s$ of $\det(\cE)$ such that
$$\widetilde{l_B}(\det(\cE),\cU,s)=\bc$$.

The notions of taming and tamed resolutions were introduced in \cite{bunke020}.
The chains $\bc$ which represent the class  $\ind^2_{Del}(\cE)$
are associated to and depend on tamed one-resolutions of $\cE$. In the following we roughly recall this construction refering to \cite{bunke020} for all details and notions. We choose a covering $\cU$ such that there exists a taming $\cE_{|U_l,t}$ for all $l\in I$. Over $U_l\cap U_k$
we have an induced  boundary taming of $\cE_{U_l\cap U_k}\times \Delta^1$.
We choose any extension to a taming $(\cE_{|U_l\cap U_k}\times\Delta^1)_t$. These choices constitute a tamed $1$-resolution $\tZ$ of $\cE$. The associated chain $\bc=\bc(\tZ) \in\check{\bZ}^1(B,\cR^1_B)$
is now given in terms of eta forms and eta invariants as follows:
\begin{eqnarray*}
c^{0,1}_l&=&\eta^1(\cE_{|U_l,t})\\
c^{1,0}_{lk}&=&\exp(-2\pi i \eta^0((\cE_{|U_l\cap U_k}\times\Delta^1)_t)) \ .
\end{eqnarray*}
 .

Now we fix a tamed one-resolution $\tZ$ and construct the local trivialization $s$. First observe that the construction of the determinant line bundle with connection and metric extend to pretamed families.
We consider the bundle $p:T\rightarrow B$ such that the fibre $T_b$ consists of
all pretamings of $\cE_b$. The model of this fiber is the space of all even selfadjoint smoothing operators, in particular, it is a Fr{\'e}chet space. Then the geometric family $\cF:= p^*\cE$ over $T$ has a tautological pretaming $\cF_t$.
Let $\theta:B\rightarrow T$ denote the zero section of $T$. For
$l\in I$ the taming $\cE_{U_l,t}$ induces another section
$T_l:U_l\rightarrow T_{|U_l}$. Let $\Gamma_l:U_l\times \Delta^1\rightarrow T$
be given by $\Gamma_l(b,t)= t T_l(b)$.
Then we obtain a pretamed family $\Gamma_l^*\cF_t$.
The restriction of $\Gamma_l^*\cF_t$ to $U_l\times \{1\}$ is tamed.
Therefore, the corresponding bundle $\det(\Gamma_l^* \cF_t)_{U_l\times \{1\}}$ is
canonically trivialized by the section $\Det(\Gamma_l^* \cF_t)$.
For $b\in U_l$ let $\gamma_b:\Delta^1\rightarrow U_l\times \Delta^1$ be the
path $\gamma_b(t)=(b,1-t)$.
Note that $\Gamma_l^* \cF_{t|U_l\times\{0\}}\cong \cE_{|U_l}$
Therefore we define the section $s_l$ of $\det(\cE_{|U_l})$ by 
$$s_l(b)=\|_{\gamma_b} \frac{\Det(\Gamma_l^* \cF_t)(b,1)}{h^{\det(\Gamma_l^*\cF)}(\Det(\Gamma_l^* \cF)(b,1))}\ ,$$
where  $\|_{\gamma_b}$ denotes the parallel transport along $\gamma_b$.
\begin{prop}
We have 
$$\widetilde l_B(\det(\cE),\cU,s)=\bc$$ 
\end{prop}
\proof
We first show that $$\frac{1}{2\pi i}\nabla^{\det(\cE)}\log s_l=\eta^1(\cE_{|U_l,t})\ .$$ Fix $b\in U_l$ and $X\in T_bB$.
Let $\sigma:(-1,1)\rightarrow U_l$ be any path with $\sigma(0)=b$ and $\sigma^\prime(0)=X$. 
Then we define
$\Sigma:(-1,1)\times\Delta^1 \rightarrow U_l\times \Delta^1$ by
$\Sigma(s,t)=(\sigma(s),t)$. The curvature
of $\Sigma^*\circ\Gamma^*_l(\det(\cF_t))$ is given by the local index form which vanishes by construction. Therefore we have
$$\nabla_X^{\det(\cE)}\log s_l=\nabla_X^{\det(\cE_{|U_l,t})} \log
 \frac{\Det(\cE_{|U_l,t})}{h^{\det(\cE_{|U_l,t})}(\Det(\cE_{|U_l,t}))}\ .$$ 
By definition of the Bismut-Freed connection the right-hand side is equal to 
$2\pi i \eta^1(\cE_{|U_l,t})$.

Next we show that $\frac{s_l}{s_k}= \exp(-2\pi i \eta^0((\cE_{|U_l\cap U_k}\times\Delta^1)_t))$. Let $b\in U_l\cap U_k$.
We can  choose a path  $P:\Delta^1\rightarrow T_b$ such that
$P(0)=T_l$, $P(1)=T_k$, and $P(s)$ is a taming of $\cE_b$ for all $s\in\Delta^1$.  We consider the map $\Sigma:\Delta^1\times\Delta^1\rightarrow T$
given by $\Sigma(s,t)=tP(s)$.
The pull-back $\Sigma^* \cF_t$ has a flat determinant bundle.
It follows, that
$$\frac{s_l}{s_k}=\exp(-2\pi i \int_{\Delta^1} P^* \eta^1(\cF_t))\ .$$
We claim that the left-hand side coincides with
$\exp(-2\pi i \eta^0((\cE_b\times\Delta^1)_t))$,
where the taming on the boundary $\cE_b\times \partial \Delta^1$ is induced by $T_l(b)$ and $T_k(b)$. To do so we first deform the taming above
to the local taming given by $P$ inside the world where both tamings and local tamings are allowed (see \cite{bunke020} for the notion of local taming). By the local variation formula we see that this does not change the $\eta^0$-invariant. Then we perform an adiabatic limit.
The function $\eta^0$ is independent of the parameter of the limit, and the result is just the integral of the $\eta^1$-form above.

This finishes the proof of the proposition.
\hB

We still fix a tamed one-resolution $\tZ$. Then we have the cycle $\bc$
and the corresponding locally trivialized bundle $(\det(\cE),\cU,s)$
such that $\widetilde l_B(\det(\cE),\cU,s)=\bc$.
The determinant section $\Det(\cE)$ can 
be represented by a family of complex valued functions $(\phi_l)_{l\in I}$ such that $\Det(\cE)_{|U_l}=\phi_l s_l$. 

\begin{prop}
We have
$$
\phi_l(b)=\left\{\begin{array}{cc} 0 & \ker(D^+_b)\not=0\\
\exp(-2\pi\imath \eta^0((\cE_b\times\Delta^1)_t))                                    \sqrt{\det(\Delta^+_b)} & \ker(D^+_b)=0
\end{array}\right.\ ,$$
where $\Delta^+_b=D_b^-D^+_b$, and where the taming
of $(\cE_b\times\Delta^1)_t$ along the boundary is induced by
$0$ and $T_l(b)$, respectively.
\end{prop}
\proof
It suffices to consider the points $b\in B$ where $\ker(D_b^+)=0$.
Since $h^{\det(\cE)}(s_l)\equiv 1$ we have
$|\phi_l|(b)=h^{\det(\cE)}(\Det(\cE))(b)=\sqrt{\det(\Delta_b^+)} $
(see \cite{berlinegetzlervergne92}, Prop. 9.41).

In order to compute the argument of $\phi_l(b)$ we deform the path
$\gamma_b:s\mapsto (1-s)T_l(b)$ to a path $\tilde \gamma_b:\Delta^1\rightarrow T_b$
with the same endpoints, but which runs through tamings.
Since $\det(\cF_t)_{|T_b}$ is flat we still have
$$s_l(b)=\|_{\tilde \gamma_b} \frac{\Det(\cF_{T_l(b)})}{h^{\det(\cF_{T_l(b)})}(\Det(\cF_{T_l(b)}))}\ .$$
We then see that
$\phi_l(b)=\exp(2\pi\imath \int_{\Delta^1} \eta^1(\tilde \gamma^* \cF_t))$.
We now again use an adiabatic limit argument in order to show that this integral is equal to $\exp(-2\pi i \eta^0((\cE_b\times\Delta^1)_t)) $.
\hB  

\section{Zero bordisms and associated sections}\label{zb}

Assume that the geometric family $\cE$ is the boundary $\partial \cW$ of a geometric family $\cW$ with boundary. We assume that $\cW$ is of product type along the boundary (see \cite{bunke020}, Sec. 2). Then it was observed by Dai-Freed \cite{daifreed94} that one can consider the exponentiated eta-invariant of $\cW$ as a section
of $\det(\cE)$. 

Fix $b\in B$.  Since $\ind(D_b^+)=0$ there exists an isometry $P:\ker D_b^+\rightarrow \ker D_b^-$. Let 
$\det(P)\in \Hom(\Lambda^{max}\ker(D^+_b),\Lambda^{max}\ker(D^-_b))$ be the induced isometry.
The choice of $P$ fixes a selfadjoint boundary condition for $D(\cW_b)$, and we let $\eta^0(\cW_b,P)$ be the corresponding eta invariant.
Then by \cite{daifreed94}, Prop. 2.15, 
$$\tau(\cW)(b):=\exp(2\pi i [\eta^0(\cW_b,P) + \frac{\dim\ker(D(\cW_b,P))}{2}])\frac{\det(P)}{h^{\det(\cE_b)}(\det(P))}\in\det(\cE)_b$$ is independent of the choice of $P$ (note the our $\eta^0$ is $1/2\eta^0$ in the usual convention).
In this way we define a unit norm section $\tau(\cW)\in\Gamma_1(\det(\cE))$.

We now fix a tamed one-resolution $\tZ$ of $\cE$ and let $\bc=\bc(\tZ)\in\check{\bZ}^1(B,\cR^1_B)$ be the corresponding cycle. Furthermore, let $s$ be the associated local trivialization of $\det(\cE)$. We now have a locally trivialized bundle with unit norm section
$(\det(\cE),\tau(\cW),\cU,s)$.

One of the main observations in \cite{bunke020} was that the zero bordism $\cW$
gives sections $\bu\in\Gamma_1(\bc)$. Let us describe the construction of $\bu=(u)$, $u\in \check{C}^0(B,\underline{U(1)}^\infty_B)$.
For each $l\in I$ we can extend the boundary taming $\cE_{|U_l,t}$ of
$\cW_{|U_l}$ to a taming $\cW_{|U_l,t}$.
Then $u$ is given by  
$$
u_l=\exp(2\pi i \eta^0(\cW_{|U_l,t}))\ .$$
Note that $u$ is independent of the choices.
We have indeed $\delta u=c^{1,0}$. 
\begin{prop}
We have $\widetilde{\Gamma l_B}(\det(\cE),\tau(\cW),\cU,s)=(\bc, \bu)\in  \Gamma \check{\bZ}^1(B,\cR^1_B)$. 
\end{prop}
\proof
From \cite{bunke020} we know that
$$\frac{1}{2\pi i} d \log u_l + c^{0,1}_l=\Omega^1(\cW_{|U_l})\ .$$
If $\phi\in\Gamma_1(\det(\cE))$ is the section which corresponds to by $\bu$, then this equation means that $\nabla^{\det(\cE)}\log \phi = 2\pi i\Omega^1(\cW)$.
On the other hand, by the generalization of \cite{daifreed94}, Thm. 1.9, to twisted Dirac operators we also have
$\nabla^{\det(\cE)}\log \tau(\cW)=2\pi i \Omega^1(\cW)$.
We conclude that the quotient $\frac{\phi}{\tau(\cW)}$ is locally constant.

We now argue that this constant must be equal to one.
Assume that $B$ is connected and that there exists $b\in B$
such that $D_b^+$ is invertible. If $b\in U_l$, then we can choose
the taming such that $T_l(b)=0$. In this case $P=0$ and
$\eta^0(\cW_{b,t})$ is equal to the eta invariant $\eta^0(\cW_b,0)$
with APS-boundary conditions. It follows that $\phi(b)=\tau(\cW)(b)$.

In the general case we can always perturb the family such that in some point $b\in B$  the operator $D_b^+$ is invertible and then apply the argument above.
\hB

Note that one can prove this proposition independently of
the knowledge of $\nabla^{\det(\cE)}\log \tau(\cW)$ in the spirit of
the proofs in Section \ref{loctr}. This would lead to an independent
verification of the computation of this derivative.

\section{Surgery}

In this section we consider two geometric families
$\cE^{\pm}$ with boundary. We assume that the geometry is of product type along the boundary.
 We assume that there are two
isomorphisms $f_i:\partial \cE^+\rightarrow \partial \cE^-$, $i=0,1$.
Then we can glue the families along the boundaries using $f_0$ and $f_1$ respectively. We denote the resulting geometric families with closed fibers by
$\cF^i:=\cE^+\cup_{f_i}(\cE^-)^{op}$.

We further form the geometric 
family
$$\cG:=(\partial \cE^+ \times \Delta^1)\cup_{f_0\cup f_1} (\partial \cE^-\times \Delta^1)^{op}\ .$$
where we glue the boundary component $\partial\cE^+\times\{i\}$ with
$(\partial\cE^-\times\{i\})^{op}$ using $f_i$.
This family is essentially the mapping torus of $f_1\circ f_0^{-1}$.

The following formula can be considered as the precise form of a surgery formula for the determinant line bundle. 
\begin{prop}\label{rr}
The line bundle
$$\det(\cF^0)\otimes \det(\cF^1)^{-1}\otimes \det(\cG)$$
admits a parallel unit-norm section $\Phi$.
\end{prop}
\proof 
We define $\Phi:=\tau(\cW)$ for a suitable zero bordism of
$\cF^0\cup_B(\cF^1)^{op}\cup_B\cG$ with $\Omega^1(\cW)=0$.
To this end we consider $S:=\Delta^1\times \Delta^1\setminus B((1/2,0),1/4)$.
We equip $S$ with the structure of a geometric manifold with corners also along the
deleted half disc.
We define 
$$\cW^{\pm}:=\cE^\pm\times \Delta^1\cup_{\sharp} (\partial \cE_{geom}^\pm \times S)^{op}
$$ 
by gluing along the common boundary
$\partial \cE^\pm\times \Delta^1 \cong  \partial \cE^\pm\times \Delta^1\times \{1\}$.
Then we glue $\cW^+$ with $(\cW^-)^{op}$
along the boundary face
$$\partial \cE^\pm\times [0,1/4]\times \{0\}\cup \partial \cE^\pm\times [3/4,1]\times \{0\}$$ using $(f_0\times \id\times\id \cup f_1\times \id\times \id )$.
The resulting family $\cW$ is the required zero-bordism.

By construction we have $\Omega^1(\cW)=0$.
 \hB

Using the section $\Phi$ in the proposition above
we can define the complex-valued function
\begin{equation}\label{e:phi}
\phi:=\frac{\Det(\cF^0)\otimes \Det(\cF^1)^{-1}\otimes \Det(\cG)}{\Phi}
\end{equation}
as a complex-valued function on $B$.

A surgery formula for the determinant section now amounts to a formula for $\phi$. The philosophy is that the determinant section is of global nature. 
In view of Proposition \ref{rr} this is in strong contrast to the determinant line bundle with its Bismut-Freed connection and Quillen metric which is of essentially local nature. A natural procedure to pin down the global nature
of the determiant sections is to perform an adiabatic limit.

In the following we investigate the adiabatic limit of
$|\phi|$. A priori it is clear that this norm is independent of the choices
involved in the construction of the bordism $\cW$. We will see on the one hand that even the adiabatic limit of $|\phi|$ depends on more than just local data coming from $\partial \cE^\pm$ and the maps $f_i$. On the other hand, the global contribution to the limit is via the scattering operator which is just an involution of a finite-dimensional vector space.

Without loss of generality we assume that $B$ consists of a single point.
We start with explaining the meaning of the adiabatic limit.
For $R>0$ we  consider the geometric manifolds
$\cE^\pm_R:=\cE^\pm\cup_\sharp \partial\cE^\pm\times [0,R]$,
where the glueing identifies $\partial \cE^\pm$ with $\partial\cE^\pm\times\{0\}$. Replacing $\cE^\pm$ by $\cE^\pm_R$ in the constructions above we obtain
$\cF^i_R$,  and $\cW_R$, respectively. By $\cG_R$ we denote the
boundary component of $\cW_R$ which corresponds to $\cG$. 
Let $\phi_R$ be the corresponding function defined by 
(\ref{e:phi}). The adiabatic limit corresponds to $R\to\infty$.
In the remainder of the present section we consider the adiabatic limit of
$|\phi_R|$. The final result is stated as Proposition 
\ref{p:adia-decom}.

Let $\Delta_0(\cE^\pm_R)$ denote the  Laplace operators $D(\cE_R^\pm)^2$
on $\cE^\pm_R$ equipped with Dirichlet boundary condition and define the integer $h:=\dim\ker(D(\partial\cE^+))$.

As an intermediate step we consider the computation \cite{PW2} of the limit   
\begin{equation}\label{e:adia}
\lim_{R\to\infty} R^{-2h} \frac{\det(\Delta(\cF^0_R))}{\det(\Delta_0(\cE^+_R))\cdot \det(\Delta_0(\cE^-_R))}= 
2^{-4h}{\det^*(\Delta(\partial\cE^+))}
\ \det(I- C_{f_0})^2\ .
\end{equation}
Here $C_{f_0}$ is a finite-dimensional operator which will be explained below,
and $\det^*$ takes the zeta regularized determinant of the operator on the complement of the kernel.
This formula holds under the condition of absense of eigenvalues of 
$\det(\Delta(\cF^0_R))$ which become exponentially small if $R$ tends to infinity. 
We are going to apply this formula also with $\cF^0_R$ replaced by $\cF^1_R$.

In order to define $C_{f_0}$ we recall the definition of the scattering matrix.
We complete $\cE^+$ by attaching an infinite cylinder :
 $\cE^+_\infty:=\cE^\pm\cup_\sharp (\partial\cE^+\times[0,\infty))$. We write the restriction of the Dirac bundle $\cV$ of $\cE_\infty^+$ to
the cylindrical part as $\cV_{|\partial\cE^+\times[0,\infty)}\cong \partial\cV * [0,\infty)$ (using the notation introduced in \cite{bunke020}).
Here the Dirac bundle $\partial\cV$  over $\partial\cE^+$ is the boundary reduction of $\cV$.
Explicitly, $\partial\cV * [0,\infty)$ is given by $\pr^*\partial{\cV}\otimes \C^2$ such that the Clifford multiplications have the form
$$c(X):=c_{\partial{\cV}}(X)\otimes \left(
\begin{array}{cc}0&i\\-i&0\end{array}\right)\  \:\:\mbox{for $X\in T\partial E^+$}\:,\:\:\:
c(\partial_u):=i\otimes \left(
\begin{array}{cc}1&0\\0&-1\end{array}\right)\ ,  
$$ and 
the grading is given by 
$$z:=1\otimes \left(
\begin{array}{cc}0&1\\1&0\end{array}\right)\ . 
$$
Here $\pr:\partial\cE^+\times[0,\infty)\rightarrow \partial\cE^+$ is the projection, and $u$ is the coordinate of $[0,\infty)$.
If $\psi$ is a section of $\partial \cV$, then we can lift it to  $u$-independent sections
$$L^{\pm}(\psi):= \psi\otimes \left(\begin{array}{c}1\\ \pm 1\end{array}\right)$$ of $\cV^\pm_{|\partial\cE^+\times[0,\infty)}$.

For $\lambda\in\R$ we now consider the space of bounded eigensections of $\Delta^\pm(\cE^+_\infty)$ to the eigenvalue $\lambda^2$. 
Let $\cH^+:=\ker(\Delta(\partial\cE^+))$ and recall that $h=\dim(\cH^+)$. 
If $\lambda^2$ is smaller than the first non-zero eigenvalue of
$\Delta^\pm(\partial\cE^+)$, then the dimension of the space of bounded
eigensections of $\Delta^\pm(\cE^+_\infty)$ is
given by $h$. All these sections have the form
$$\ee^{-i\lambda u} L^\pm (\psi) +\ee^{i\lambda u} L^\pm(C^\pm_{\cE^+}(\lambda)\psi)+o(u)$$
with uniquely determined $\psi\in \cH^+$, and
for a certain operator $C^\pm_{\cE^+}(\lambda)\in\End(\cH^+)$. This operator is called the scattering operator. Using the fact that
$D(\cE^+)$ maps eigensections of $\Delta^+(\cE^+)$ to eigensections of
$\Delta^-(\cE^+)$ one can check 
that $C^+_{\cE^+}(\lambda)=-C^-_{\cE^+}(\lambda)$.
Furthermore, the scattering operator satisfies the functional equation
 $$C^\pm_{\cE^+}(\lambda)C^\pm_{\cE^+}(-\lambda)=\id\ ,$$
which implies that $C^\pm_{\cE^+}(0)$ is an involution.
Note that $f_0:\partial\cE^+\rightarrow\partial\cE^-$ induces an isomorphism $\cH(f_0):\cH^+\rightarrow\cH^-$.
Now we can define
$$C_{f_0}:= \cH(f_0)^{-1}\circ C^-_{\cE^-}(0)\circ \cH(f_0)\circ C^+_{\cE^+}(0)\in \End(\cH^+)\ .$$   
In a similar manner we define
$C_{f_1}\in\End(\cH^+)$.
Note that $C_{f_i}^*= C^+_{\cE^+}(0)\circ C_{f_i}\circ C^+_{\cE^+}(0)$.
We conclude that the non-real eigenvalues of the unitary operator
$C_{f_i}^*$ come in pairs $(\mu,\bar \mu)$ with $|\mu|=1$, $\Imm(\mu)>0$.
Because of the assumption about the absence of exponentially decreasing  eigenvalues $1$ is not in the spectrum of $C_{f_i}$.
We conclude that $\det(1-C_{f_i})>0$.

Using that $$h^{\det(\cF^i_R)}(\Det(\cF^i_R))=\det(\Delta^+(\cF^i_R))=\sqrt{\det(\Delta(\cF^i_R))}$$ we conclude
that $$\lim_{R\to\infty} \frac{h^{\det(\cF^0_R)}(\Det(\cF^0_R))}{h^{\det(\cF^1_R)}(\Det(\cF^1_R))}=
\frac{\det(1-C_{f_0})}{\det(1- C_{f_1})}\ .$$

We now consider $\det(\Delta(\cG_R))$ in the adiabatic limit. Note that
$$\cG_R:=(\partial \cE^+ \times \Delta_R^1)\cup_{f_0\sqcup f_1} (\partial \cE^-\times \Delta_R^1)^{op}\ ,$$
where $\Delta^1_R\cong [0,R]$. For the adiabatic limit of the determinant of $\det(\Delta(\cG_R))$, we have
the following formula in \cite{PW3},
\begin{equation}\label{e:adia2}
 \det(\Delta(\cG_R))\stackrel{R\to\infty}{\sim} \det(\Delta_{f_0\sqcup f_1})^2\cdot \exp\left(\ 4R \left[\ (\psi(-\frac {1}{2})+\gamma)\ a_{\frac{1}{2}}
+\zeta^r_{\Delta(\partial\cE^+)}(-\frac 12)\right]\ \right)\ . 
\end{equation}
Here the notation $\stackrel{R\to\infty}{\sim} $ means that the quotient of each side converges to $1$ as $R\to\infty$.
Note that the second term may blow up in general as $R\to\infty$. 
Here $\Delta_{f_0\sqcup f_1}$ is the Laplacian acting on the flat vector bundle with the holonomy $\cH(f_0)^{-1}\circ \cH(f_1)$
over $S^1$. For the second term,
$\psi(x)=\frac{d}{dx} \log \Gamma(x)$, $\gamma$ is  Euler's constant 
and $a_{\frac{1}{2}}$ and $\zeta^r_{\Delta(\partial\cE^+)}(s)$  are the residue and regular parts 
of $\zeta_{\Delta(\partial\cE^+)}(s)$ at $s=-\frac 12$ such that
$$
\zeta_{\Delta(\partial\cE^+)}(s)\ =\ \frac{a_{\frac12}}{s+\frac 12}  + \ \zeta^r_{\Delta(\partial\cE^+)}(s) \ .
$$

We eventually conclude 

\begin{prop}\label{p:adia-decom}
If $\Delta(\cF_R^i)$ and $\Delta(\cG_R)$ have no exponentially decreasing eigenvalues, then we have the following asymptotic behaviour of $|\phi_R|$ in the adiabatic limit :
$$|\phi_{R}| \stackrel{R\to\infty}{\sim} \ 
\det(\Delta_{f_0\sqcup f_1})\frac{\det(1-C_{f_0})}{\det(1-C_{f_1})}\exp( 2R\ [\ (\psi(-\frac {1}{2})+\gamma)\ a_{\frac{1}{2}}
+\zeta^r_{\Delta(\partial\cE^+)}(-\frac 12)\ ]\ )
\ .$$
\end{prop}

Note that the right hand side of this formula contains the term
$\frac{\det(1-C_{f_0})}{\det(1-C_{f_1})}$ which in general depends on the global structure of $\cE^\pm$.

It might be an interesting problem to study the adiabatic limit of the phase
of $\phi$.

{\bf Remark :} The surgery formula for the determinant line bundle with Quillen metric and Bismut-Freed connection can also be deduced from the splitting formula of Piazza \cite{piazza98}. On the other hand, Proposition
\ref{p:adia-decom} does not follow from this work.

In \cite{scott00}, Scott also proved a splitting formula for the determinant line  bundle. But he considers a different metric and connection which are associated to Fredholm determinants rather than zeta regularized determinants.

Assume that we have fixed an isomorphism $\partial \cE^+\cong \partial\cE^-$.
The idea to associate a determinant line bundle to the families $\cE^\pm$  and
to express the determinant line bundle of $\cE=\cE^+\cup_\sharp (\cE^-)^{op}$ in terms of these pieces (this is the content of a splitting formula) in order to prove a surgery formula seems to be not so natural.
The point is that in order to define the determinant line bundle
for $\cE^\pm$ one must choose  non-natural boundary conditions. In the remainder of the
present paper we want to advertise the idea that $\cE^\pm$ give rise
to trivializations of the index gerbes of $\partial\cE^\pm$.
The spitting formula for the determinant line bundle is then equivalent
to the statement that the line bundle (with metric and connection) given by the difference of the two trivializations
of the index gerbe of $\partial\cE^+\cong\partial \cE^-$ induced by
$\cE^+$  and $\cE^-$  is isomorphic to the determinant line bundle of $\cE$.
{\em Note that we talk here about the determinant line bundle with the Quillen metric and the Bismut-Freed connection and not about the determinant sections .}

\section{Gerbes and sections}\label{bound}

The goal of the present and the following section is to explain that the generalization
of the determinant line bundle to families  with boundaries leads to
a new kind of object : sections of gerbes.
The starting point is the observation that there is no canonical choice of
boundary conditions in general. The way out is to consider all choices
(in some class) at once.

In the preceding sections we already met a similar problem in one degree less.
Here the exponentiated eta-invariant of a manifold depends on the choice of boundary conditions. In order to make it natural, we were led to consider it as a section of the determinant line bundle associated to the boundary.

We consider the complex of sheaves
$$\cR^2_B:  
 0\rightarrow \underline{U(1)}^\infty_B\stackrel{\frac{1}{2\pi i}d \log }{\rightarrow}\cA^1_B\rightarrow \cA_B^2\rightarrow 0\ .$$
Then a gerbe is by definition a cycle
$\bc\in \check{\bZ}^2(B,\cR^2_B)$. It is a tuple
$\bc=(c^{0,2},c^{1,1},c^{2,0})$, where
\begin{eqnarray*}
c^{0,2}&\in&\check{C}^0(B,\cA^2_B)\\
c^{1,1}&\in&\check{C}^1(B,\cA^1_B)\\
c^{2,0}&\in&\check{C}^2(B,\underline{U(1)}_B^\infty)\ .
\end{eqnarray*}
The cycle condition reads
\begin{eqnarray*}
\delta c^{2,0}&=&0\\
\delta c^{1,1}&=&\frac{1}{2\pi i} d \log c^{2,0}\\
\delta c^{0,2}&=&- d c^{1,1}\ .
\end{eqnarray*}
An isomorphism $\bb:\bc_0\rightarrow \bc_1$
is given by a chain $\bb=(b^{0,1},b^{1,0})\in \check{\bC}^1(B,\cR^2_B)$, where
\begin{eqnarray*}
b^{0,1}&\in&\check{C}^0(B,\cA^1_B)\\
b^{1,0}&\in&\check{C}^1(B,\underline{U(1)}_B^\infty)
\end{eqnarray*}
such that
\begin{eqnarray*}
d b^{0,1}&=&c^{0,2} \\
- \frac{1}{2\pi i}d \log b^{1,0} - \delta b^{0,1}&=& c^{1,1}\\
-\delta b^{1,0}&=& c^{2,0}\ .
\end{eqnarray*}
By definition, $H^3_{Del}(B)$ is the group of isomorphism classes of gerbes.

The curvature of the gerbe is the three-form $R^\bc\in\cA_B^3(B,d=0)$ which characterized by $R^{\bc}_{|U_l}=d c^{0,2}_l$. Furthermore, there is  a characteristic class $\bv(\bc)\in \check{H}^3(B,\underline{\Z}_B)$, which is represented by $\delta \frac{1}{2\pi i} \log c^{2,0}$.
Note that both, the curvature and the characteristic class, only depend on the isomorphism class of the gerbe. A section of a gerbe is a chain
$\bb=(b^{0,1},b^{1,0})$ as above satisfying only
\begin{eqnarray*}
- \frac{1}{2\pi i}d \log b^{1,0} - \delta b^{0,1}&=& c^{1,1}\\
-\delta b^{1,0}&=& c^{2,0}\ .
\end{eqnarray*}
Let $\Gamma_1(\bc)$ denote the set of all sections. 
Note that an isomorphism $\bb:\bc_0\rightarrow \bc_1$ of gerbes induces an isomorphism
$\Gamma(\bb):\Gamma_1(\bc_0)\rightarrow\Gamma_1(\bc_1)$ of sections by
$\Gamma(\bb)(\bb_0):=\bb_0+\bb$.
We let $\Gamma\check{\bZ}^2(B,\cR^2_B)$ be the set of pairs $(\bc,\bb)$, where $\bb\in\Gamma_1(\bc)$. Further, let $\Gamma H^3_{Del}(B)$ denote the set of isomorphism classes in $\Gamma\check{\bZ}^2(B,\cR^2_B)$.

If $\bb\in \Gamma_1(\bc)$,
then there is a two-form $\nabla^\bc\bb\in\cA^2(B)$ such that
$$\nabla^\bc\bb\in\cA^2(B)_{|U_l}= d b^{0,1}_l - c^{0,2}_l\ .$$
We call two sections $\bb_i$, $i=0,1$, isomorphic if there is a chain
$u\in \check{C^0}(B,\underline{U(1)}_B^\infty)$ such that
\begin{eqnarray*}
\delta u&=&b^{1,0}_1-b^{1,0}_0\\
\frac{1}{2\pi i}d \log u&=&b^{0,1}_1-b^{0,1}_0\ .
\end{eqnarray*}

Observe, that the difference of two sections $\bb_i\in\Gamma_1(\bc)$ is a an element of $\check{\bZ}^1(B,\cR^1_B)$, i.e. it corresponds to an isomorphism class
of locally trivialized line bundles. The curvature of this line bundle is given by $\nabla^\bc(\bb_1-\bb_0)$. Similarly, the difference of two isomorphism classes of sections of the gerbe $\bc$ is just a class in $H^2_{Del}(B)$.
We thus have a pairing
$$<.,.>:\Gamma \check{\bZ}^2(B,\cR^2_B)\times_{ \check{\bZ}^2(B,\cR^2_B)}\Gamma  \check{\bZ}^2(B,\cR^2_B)\rightarrow H^2_{Del}(B)$$ given by $<(\bc,\bb_0),(\bc,\bb_1)>:=[\bb_1-\bb_0]$.
Note that this pairing does not factor over the isomorphism classes $[\bc_i,\bb_i]$.

\section{The determinant line bundle of a family with boundary}

We  consider an even-dimensional geometric family $\cE$ with boundary $\partial \cE$. 
Again, we assume that $\cE$ is of product type along the boundary. 
$\partial \cE$ is a geometric family with closed odd-dimensional
fibers. We have $\ind(\cE)=0$ in $K^1(B)$. The index gerbe of such a family was first introduced in \cite{lott01}. In the present note we adopt the point of view of \cite{bunke020}, where the isomorphism class of the index gerbe is given by $\ind^3_{Del}(\partial \cE)\in H^3_{Del}(B)$.
A gerbe $\check{\bZ}^2(B,\cR^2_B)$ representing this class is associated to the
choice of a tamed two-resolution $\tZ$.
This resolution is given by the following choices. First we need
an open cover $\cU=(U_l)_{l\in I}$ and a taming $\partial\cE_{|U_l,t}$ for all $l\in I$. Then we need an extension of the induced boundary taming of $\partial \cE_{|U_l\cap U_k}\times \Delta^1$
to a taming $(\partial \cE_{|U_l\cap U_k}\times \Delta^1)_t$ for all $l,k\in I$ with $U_l\cap U_k\not=\emptyset$.
Finally, we need an extension of induced boundary taming
of $\partial \cE_{|U_l\cap U_k\cap U_h }\times \Delta^2$
to a taming $(\partial \cE_{|U_l\cap U_k\cap U_h }\times \Delta^2)_t$ for all $l,k,h\in I$ with $U_l\cap U_k\cap U_h\not=\emptyset$.
Then $\bc=\bc(\tZ)$ is given by (see \cite{bunke020})
\begin{eqnarray*}
c^{0,2}_l&:=&\eta^2(\partial\cE_{|U_l,t})\\
c^{1,1}_{lk}&:=&-\eta^1((\partial \cE_{|U_l\cap U_k}\times \Delta^1)_t)\\
c^{2,0}_{lkh}&:=&\exp(2\pi i \eta^0((\partial \cE_{|U_l\cap U_k\cap U_h }\times \Delta^2)_t))\ .
\end{eqnarray*}
We now consider the family $\cE$ as a zero-bordism of $\partial \cE$. It extends to a zero bordism $\tW$ of the tamed two-resolution $\tZ$.
In order to ensure that such a lift exists we must choose the taming
of $\partial \cE_{|U_l,t}$ such that $\cE_{|U_l,bt}$
has trivial index for all $l$. Such a choice is possible.
$\tW$ is given by extension of the boundary tamings to tamings $\cE_{|U_l,t}$, and further, extensions of the boundary tamings of
$\cE_{|U_l\cap U_k}\times\Delta^1$ to tamings $\cE_{|U_l\cap U_k}\times\Delta^1$. Then we have the section
$\bb:=\bb(\tW)\in \Gamma_1(\bc)$ given by
\begin{eqnarray*}
b^{0,1}_l&:=&-\eta^1(\cE_{|U_l,t})\\
b^{1,0}_{lk}&:=&\exp(-2\pi i \eta^0((\cE_{|U_l\cap U_k}\times \Delta^1)_t))\ .
\end{eqnarray*}
Note that $\nabla^\bc \bb=-\Omega^2(\cE)$.
The following proposition generalizes from degree two to degree three 
the observation that a geometric family with boundary gives rise to a well-defined isomorphism class of a locally trivialized line bundle with section, namely, the determinant line bundle of the boundary with the section introduced in Section \ref{zb}.
\begin{prop}
The pair  $(\bc,\bb)\in \Gamma\check{\bZ}^2(B,\cR^2_B)$ is independent of all choices of tamings.
\end{prop}
\proof
Assume that we have two tamed two-resolutions $\tZ_0$ and $\tZ_1$ of $\partial \cE$ given without loss of generality with respect to the same covering $\cU$ with index set $I$.
Then we consider the covering $\cV$ with index set $J:=I\times\{0,1\}$
where $V_{(l,\epsilon)}:=U_l$.
We define a tamed two resolution $\tZ^\prime$ with respect to $\cV$ as follows. On $\partial \cE_{|U_{l,\epsilon}}$ we choose the taming
given by $\tZ_\epsilon$. In the next step we choose
$(\partial \cE_{|U_{l,\epsilon}\cap U_{k,\delta}}\times \Delta^1)_t$
as given by $\tZ_\epsilon$, if $\epsilon=\delta$, and arbitrary else.
This is again possible since we require trivial index
of $\cE_{|U_{l,\epsilon},bt}$. We continue to define $\tZ^\prime$ on triple intersections in the same manner.

Now the main point is that $\tZ_\epsilon$, $\epsilon\in\{0,1\}$, are both obtained by a refinement
of $\tZ^\prime$ so that $\bc(\tZ_0)=\bc(\tZ_1)\in \check{\bZ}^2(B,\cR^2_B)$.
In the same manner given the zero bordisms $\tW_\epsilon$ of $\tZ_\epsilon$ we can find a zero bordism $\tW^\prime$ of $\tZ^\prime$ with respect to $\cV$ such that the two zero bordisms $\tW_\epsilon$ are obtained by refinements. Then we have
$\bb(\tW_0)=\bb(\tW_1)\in \check{\bC}^1(B,\cR^2_B)$. 
\hB

Now assume that we have two geometric families $\cE^\pm$ with boundary
$\partial \cE^\pm$ as above together with an isomorphism
$\partial \cE^+\cong \partial \cE^-$. Then we form the glued family
$\cE:=\cE^+\cup_\sharp(\cE^-)^{op}$. We  have two pairs
$\gamma(\cE^\pm,\partial \cE^\pm)\in \check{\bZ}^2(B,\cR^2_B) $.
Let us assume for simplicity that the numerical index of $\cE$ vanishes.
The gluing formula for the determinant bundle reads
\begin{prop}
In $H^2_{Del}(B)$ we have the identity 
$$<\gamma(\cE^+,\partial \cE^+),\gamma(\cE^-,\partial \cE^-)>=\det(\cE)\ .$$
\end{prop}
\proof
We construct a suitable bordism $\cW$.
Let $S\subset \R^2$ be a pentagon. Let $\partial_i$, $i=1,\dots 5$, be the boundary faces in cyclic order.
We equip $S$ with the structure of a geometric manifold with corners such that
all boundary components $\partial_i S$ are isometric to $I:=[0,1]$.
We form the geometric family $\cW$ by gluing
$(\cE^+\times I) \cup_\sharp \partial \cE^+\times S\cup_\sharp (\cE^-\times I)$ according to  $(\partial \cE^+\times I)\stackrel{\sim}{\rightarrow} \partial_1 (\partial \cE^+\times S)$, 
$\partial_4 (\partial \cE^+\times S)\stackrel{\sim}{\rightarrow}(\partial \cE^-\times I)^{op}$.
The resulting geometric family $\cW$ has three boundary faces which are isomorphic to $\cE^\pm$ and $\cE$.

We fix a tamed two-resolution $\tZ^+$  of $\partial\cE^+$.
It induces a tamed two-resolution $\tZ^-$ of $\partial \cE^-$.
Then we fix tamed  zero-bordisms $\tW^\pm$ of $\tZ^\pm$ associated to
the zero-bordisms $\cE^\pm$ of $\partial\cE^\pm$.
We also fix a tamed one-resolution $\tU$ of $\cE$. This induces a boundary taming of $\cW_{|U_l}$ which we extend to a taming $\cW_{|U_l,t}$ for each $l\in I$.

Let $\bc:=\bc(\tZ)$ and $\bb^\pm:=\bb(\tW^\pm)\in\Gamma_1(\bc)$.
Furthermore, let $\bb:=\bb(\tU)\in \check{\bZ}^1(B,\cR^1_B)$ as in Section \ref{loctr}. We must construct an isomorphism
$u:\bb_1-\bb_0\rightarrow \bb$.
We will consider $u\in\Gamma_1(\bb+\bb_0-\bb_1)$ such that
$\nabla^{\bb+\bb_0-\bb_1}u=0$.
Let $u$ be given by 
$$u_l:=\exp(-2\pi i \eta^0(\cW_{|U_l,t}))\ .$$
Then indeed
$(\delta u)_{lk}=[\bb+\bb_0-\bb_1]^{1,0}_{lk}$ and
$d u_l = [\bb+\bb_0-\bb_1]^{0,1}_{l}$ as required.
\hB

\bibliographystyle{plain}

\begin{thebibliography}{10}

 
 
\bibitem{berlinegetzlervergne92}
N.~Berline, E.~Getzler, and M.~Vergne.
\newblock {\em Heat Kernels and Dirac Operators}.
\newblock Springer-Verlag Berlin Heidelberg New York, 1992.


\bibitem{bismutfreed861}
J.~M. Bismut and D.~Freed.
\newblock The analysis of elliptic families. I. Metrics and connections on
  determinant bundles.
\newblock {\em Comm. Math. Phys.}, 106(1986), 159--176.

 
\bibitem{brylinski93}
J.~L. Brylinski.
\newblock {\em Loop spaces, characteristic classes, and geometric
  quantization}.
\newblock Birk{\"a}user, Progress in Math. 107, 1993.

\bibitem{bunke020}
U.~Bunke.
\newblock Index theory, eta forms, and Deligne cohomology.
\newblock arXiv:math.DG/0201112

\bibitem{daifreed94}
X.~Dai and D.~Freed.
\newblock $\eta$-invariants and determinant lines. Topology and physics.
\newblock {\em J. Math. Phys.}  35(1994), 5155--5194. arXiv:hep-th/9405012


\bibitem{gajer97}
P.~Gajer.
\newblock Geometry of Deligne cohomology.
\newblock {\em Invent. Math.}, 127(1997), 155--207. arXiv:alg-geom/960125

 
\bibitem{hitchin99}
N.~Hitchin.
\newblock Lectures on special Lagrangian submanifolds.
\newblock arXiv:math.DG/9907034

\bibitem{leinster01}
T.~Leinster.
\newblock A Survey of Definitions of n-Category.
\newblock arXiv:math.CT/0107188


\bibitem{lott01}
J.~Lott.
\newblock Higher degree analogs of the determinant line bundle.
\newblock {\em Comm. Math. Phys.}, 230(2002), 41-69. arXiv:math.DG/0106177

\bibitem{murray96}
J. Murray. 
\newblock Bundle gerbes.  
\newblock{\em J. London Math. Soc.(2)}  54 (1996),  403--416.
arXiv:math.dg-da/9407015 


\bibitem{PW2}
J.~Park and K. ~Wojciechowski.
\newblock Scattering theory, the adiabatic decomposition of the $\zeta$-determinant
and the Dirichlet to Neumann operator.
\newblock preprint

\bibitem{PW3}
J.~Park and K.~Wojciechowski.
\newblock Holonomy theorems for $\eta$-invariant and $\zeta$-determinant.
\newblock preprint


\bibitem{piazza98}
P.~Piazza.
\newblock Determinant bundles, manifolds with boundary and surgery. II. Spectral sections and surgery rules for anomalies.
\newblock {\em Comm. Math. Phys}. 193(1998), 105--124. 



\bibitem{scott00}
S.~Scott.
\newblock Splitting the curvature of the determinant line bundle.
\newblock {\em Proc. Amer. Math. Soc}. 128 (2000), 2763--2775. arXiv:math.AP/9812124

\end{thebibliography}

\end{document}